\documentclass[12pt]{amsart}
\usepackage{amsfonts}
\usepackage{amssymb}
\theoremstyle{definition}

\theoremstyle{remark}

\begin{document}

\vspace*{0.75 in}

\centerline{\bf Comptes rendus de l'Acad\'emie bulgare des Sciences}

\centerline{\it Tome 54, No 3, 2001}

\vspace{0.6in}

\begin{flushright}
{\it MATH\'EMATIQUES
\\ G\'eometrie diff\'eretielle}
\end{flushright}

\vspace{0.2in}

\par\vspace {1cm}
\centerline{\large \bf   SCHUR'S THEOREM FOR ALMOST      }
\centerline{\large \bf   HERMITIAN MANIFOLDS  }
\par\vspace{5 ex}
\centerline{\bf O. T. Kassabov}
\par\vspace {4 ex}

\vspace{0.2in}
\centerline{\it (Submited by Academician P. Kenderov on October 11, 2000)}

\vspace{0.2in}

\par
\section {\bf Introduction.}
Let $M$ be a Riemannian manifold with curvature tensor $ R $. The sectional
curvature of a 2-plane $ \alpha $ in a tangent space $T_pM $ is defined
by
$$ K(\alpha ,p) = R(x,y,y,x) \ \ , $$
where $ \{ x,y \} $ is an orthonormal basis of $ T_pM $. The classical
theorem of F. Schur says that if $ M $ is a connected manifold of dimension
$ n \ge 3 $ and in any point $ p \in M $ the curvature $ K(\alpha ,p ) $
does not depend on $ \alpha \in T_pM $, then it does not depend on the
point $ p $ too, i.e. it is a global constant. Such a manifold is called a
manifold of constant sectional curvature.
\par Now let $ M $ be an almost Hermitian manifold with metric tensor $ g
$, almost complex structure $ J $ and curvature tensor $ R $. A 2-plane
$ \alpha \in T_pM $ is called holomorphic, resp. antiholomorphic if
$ J\alpha = \alpha $, resp. $ J\alpha \perp \alpha $. Then $ M $ is said to
be of pointwise constant holomorphic (resp. antiholomorphic)
sectional curvature $ \nu (p) $, if in any point $ p\in M $ $ K(\alpha ,p)
= \nu (p) $ for any holomorphic (resp. antiholomorphic) 2-plane $ \alpha
\in T_pM $. Otherwise, in any point of $ M $ the curvature of any
holomorphic (resp. antiholomorphic) 2-plane does not depend on $ \alpha $.
\par For the special class of K\"ahler manifolds both the holomorphic and
the antiholomorphic analogues of the Schur's theorem hold, see [1], [4].
In fact, for K\"ahler manifolds the requirements of pointwise constant
holomorphic and pointwise constant antiholomorphic sectional curvature are
equivalent [1].
\par In [3] we proved the Schur's theorem of antiholomorphic type for the
so-called $ RK $-manifolds or $ AH_3 $-manifolds, i.e. almost Hermitian
manifolds, whose curvature tensor satisfies
$$ R(x,y,z,u) = R(Jx,Jy,Jz,Ju)  \ \ . $$
In this note we shall prove a Schur's theorem for an arbitrary almost
Hermitian manifold. Namely, we have:
\par {\bf Theorem.} {\it Let $ M $ be a connected $2n$-dimensional almost
Hermitian manifold, $ n \ge 3 $. If $ M $ is of pointwise constant
antiholomorphic sectional curvature $ \nu (p) $, then $ \nu $ is a global
constant.}

\vspace{0.2in}
\section {\bf Preliminaries. }
For a tensor field $ Q $ of type (0,2) denote by $ \varphi (Q) $ and
$ \psi (Q) $ the tensor fields of type (0,4) defined respectively by
$$ \begin{array} {rl}
\varphi (Q)(x,y,z,u) & = g(x,u)Q(y,z) - g(x,z)Q(y,u) \\
                     & + g(y,z)Q(x,u) - g(y,u)Q(x,z)
\end{array}  \ \ \ , $$
$$ \begin{array} {rl}
\psi (Q)(x,y,z,u) & = g(x,Ju)Q(y,Jz)-g(x,Jz)Q(y,Ju)-2g(x,Jy)Q(z,Ju) \\
                 & +  g(y,Jz)Q(x,Ju)-g(y,Ju)Q(x,Jz)-2g(z,Ju)Q(x,Jy)
\end{array}  \ \ \ . $$
Denote also $ \pi_1 = \frac{1}{2} \varphi (g) $,
$ \pi_2 = \frac{1}{2} \psi (g) $.
\par Let $ \{ e_1,e_2,...,e_{2n} \} $ be an orthonormal basis of
$ T_pM $. The Ricci $ * $-tensor $ \rho^* $ is given by
$$ \rho^* (x,y) = \sum_{i=1}^{2n} R(x,e_i,Je_i,Jy) $$
and the $ * $-scalar curvature $ \tau^* $ is
$$ \tau^* (p) = \sum_{i=1}^{2n} \rho^*(e_i,e_i) \ \ .  $$
The curvature tensor of a $2n$-dimensional almost Hermitian manifold of
pointwise constant antiholomorphic sectional curvature $ \nu (p) $ has the
form [2]
$$ R = \frac{1}{2(n+1)} \psi (\rho^*) + \nu \pi_1 -
       \frac{\tau^* +2(n+1)\nu}{2(n+1)(2n+1)} \pi_2 \ \ . $$
We can rewrite this in the form
$$ R =  \psi (Q) + \nu \pi_1   \ \ ,  \leqno (2.1)  $$
where
$$ Q = \frac{1}{2(n+1)} \rho^* -
     \frac{\tau^* +2(n+1)\nu}{4(n+1)(2n+1)}g  \ \ . $$
Note, that the tensor $ Q $ is neither symmetric, nor $ J $-invariant, but
it has the property
$$  Q(Jx,Jy) = Q(y,x) \ \ . $$

\vspace{0.2in}
\section {\bf Proof of the Theorem. }
Denote by $ s(Q) $ the symmetric part of $Q$. In a point $ p \in M $ we
choose an orthonormal basis $ \{ e_1,Je_1,...,e_n,Je_n \} $ of $ T_pM $
such that
$$ s(Q)(e_i,e_j) = s(Q)(e_i,Je_j) = 0 \ \ \ {\rm for}\ \ i \ne j \ . $$
Then $ Q(e_i,Je_j) = -Q(e_j,Je_i) $. In what follows $ (x,y,z) $ will stand
for an antiholomorphic triple of vectors, bellonging to the above chosen
basis.

\vspace{0.1in}
\par {\bf Lemma 1.} {\it There exists a non-zero vector
$ \tilde y \in span\{ y,Jy \} $, such that $ Q(x,J\tilde y) = 0 $. }

\vspace{0.1in}
\par Proof. Suppose $ Q(x,Jy) \ne 0 $. Then the vector
$$ \tilde y = \frac{Q(x,y)}{Q(x,Jy)} y + Jy $$
satisfies the desired condition.

\vspace{0.1in}
\par {\bf Lemma 2.} {\it There exist non-zero vectors
$ \tilde x \in span\{ x,Jx \} $, $ \tilde y \in span\{ y,Jy \} $,
$ \tilde z \in span\{ z,Jz \} $, such that }
$$ Q(\tilde x,J\tilde y) = 0 \qquad
   Q(\tilde x,J\tilde z) = 0 \qquad
   Q(\tilde y,J\tilde z) = 0    \ \ .   \leqno(3.1)$$

\vspace{0.1in}   
\par Proof. Using Lemma 1 we can assume that $ Q(x,Jy) = Q(x,Jz) = 0 $.
Suppose $ Q(y,Jz) \ne 0 $. Let $ \tilde x = ax + Jx $, $ \tilde y = by + Jy
$, $\tilde z = cz + Jz $. Then the conditions (3.1) are satisfied if
$$ \left\{ \begin{array}{l}
a+b=0  \\
a+c=0  \\
(b+c)Q(y,z) +(1-bc)Q(y,Jz) = 0 \ \ .
\end{array}  \right.  $$
The last system having always a solution, the Lemma is proved.

\vspace{0.1in}
\par Now we prove the Theorem. Denote by $ \nabla $ the covariant
differentiation on $ M $. From the second Bianchi identity
$$ (\nabla_xR)(y,z,z,y) + (\nabla_yR)(z,x,z,y)
                        + (\nabla_zR)(x,y,z,y) = 0 \ \ , $$
using (2.1) we obtain
$$  \begin{array} {rl}
  x(\nu ) - 6g((\nabla_xJ)y,z)Q(y,Jz)  &   \\
 + 3\{ g((\nabla_yJ)y,z)Q(x,Jz)  + g((\nabla_yJ)x,z)Q(y,Jz)  &   \\
 + g((\nabla_zJ)z,y)Q(x,Jy) + g((\nabla_zJ)y,x)Q(y,Jz) \} & =0 \ .  \\
\end{array}   \leqno (3.2)   $$
Using Lemma 2 we may assume that $ Q(x,Jy)=Q(x,Jz)=Q(y,Jz) = 0 $.
Then (3.2) implies $ x(\nu ) = 0 $. Now we change $ (x,y,z) $ in (3.2) by
$ (Jx,Jy,Jz) $ to find
$$  \begin{array} {rl}
  Jx(\nu ) - 6g((\nabla_{Jx}J)y,z)Q(y,Jz)  &   \\
 + 3\{ g((\nabla_{Jy}J)y,z)Q(x,Jz)  + g((\nabla_{Jy}J)x,z)Q(y,Jz)  &   \\
 + g((\nabla_{Jz}J)z,y)Q(x,Jy) + g((\nabla_{Jz}J)y,x)Q(y,Jz) \} & =0   \\
\end{array}     $$
and hence $ Jx (\nu ) = 0 $. Consequently $ \nu $ does not depend on the
vectors in $ span \{x,Jx \} $ and thus it is a global constant.

\vspace{1 cm}
\centerline { REFERENCES }
\vspace{1 cm}
\par 1. B.-Y. Chen, K. Ogiue. Duke Math. J. {\bf 40},(1973), 797-799.
\par 2. G. Ganchev. Pliska {\bf 9}(1987), 33-43.
\par 3. O. T. Kassabov. C. R. Acad. bulg. Sci. {\bf 35},(1982), 905-907.
\par 4. S. Kobayashi, K. Nomizu. Foundations of Differential Geometry II,
\par \ \ \ \ New York, 1969.

\vspace {1cm}
\par Higher Transport School "T. Kableshkov"
\par Section of Mathematics
\par Slatina, 1574 Sofia, BULGARIA

\end{document}